\documentclass[
  onefignum,onetabnum]{siamart171218}

\usepackage{amstext}
\usepackage{stmaryrd}
\usepackage{algpseudocode}

\usepackage{units}
\usepackage{cancel}
\usepackage{graphicx}
\usepackage{subcaption}

\usepackage{xcolor}
\usepackage{array}
\usepackage{verbatim}
\usepackage{booktabs}
\usepackage{calc}
\usepackage{textcomp}
\usepackage{multirow}
\usepackage{tikz,tikzscale}

\usepackage{pgfplots}
\usepackage{pgfkeys}
\pgfplotsset{compat=1.18}

\definecolor{deepgreen}{RGB}{0,100,0}

\newcommand{\plotwidth}{0.49\columnwidth}
\newcommand{\plotheight}{0.4\columnwidth}

\pgfplotscreateplotcyclelist{paper markers}{
  {teal,    solid, mark=*},
  {orange,   solid, mark=triangle*},
  {blue,     solid, mark=square*},
  {red,      solid, mark=pentagon*},
  {violet,     solid, mark=otimes*},
  {brown,    solid, mark=halfcircle*},
  {black,   solid, mark=diamond*},
}

\pgfplotsset{
  paperplot/.style={
      width=\plotwidth,
      height=\plotheight,
      cycle list name=paper markers,
      every axis/.append style={font=\small},
      grid=major,
      thick
    }
}

\pgfplotsset{
  longplot/.style={
      width=0.8\columnwidth,
      height=\plotheight,
      cycle list name=paper markers,
      every axis/.append style={font=\small},
      grid=major,
      thick
    }
}

\usepackage[all]{xy}

\usepackage{mathtools}
\usepackage{placeins}

\usetikzlibrary{calc}
\usetikzlibrary{shadings}

\usepgfplotslibrary{fillbetween}

\usepackage[mode=multiuser,status=draft]{fixme} 
\fxsetup{innerlayout=inline}
\fxsetup{targetlayout=colorcb}

\fxusetheme{color}
\FXRegisterAuthor{mw}{envmw}{MW}

\definecolor{apply}{rgb}{0.3,.7,0.}
\definecolor{applyface}{rgb}{.9,.95,0.}
\definecolor{invert}{rgb} {0.5,0.,1.}
\colorlet{invertface}{invert!60!red}

\usepackage[colorinlistoftodos,prependcaption,textsize=small]{todonotes}

\def\mesh{\mathbb M}

\def\R{\mathbb R}

\def\resth#1{I_{#1}^{\downarrow_{h}}}
\def\prolh#1{I_{#1}^{\uparrow_{h}}}

\def\prolp#1{I_{#1}^{\Uparrow_{p}}}
\def\restp#1{I_{#1}^{\Downarrow_{p}}}

\usepackage{colortbl} 
\newcommand{\myrowcolor}{\rowcolor[gray]{0.925}}

\usepackage{booktabs} 
\usepackage[font=footnotesize,position=b]{subcaption}
\usepackage{adjustbox}

\newcommand{\pluseq}{\stackrel{+}{\gets}}

\newcommand{\LUpNew}{local\_update}

\makeatletter
\def\pgfplots@getautoplotspec into#1{%
    \begingroup
    \let#1=\pgfutil@empty
    \pgfkeysgetvalue{/pgfplots/cycle multi list/@dim}\pgfplots@cycle@dim
    \let\pgfplots@listindex=\pgfplots@numplots
    \pgfkeysgetvalue{/pgfplots/cycle list set}\pgfplots@listindex@set
    \ifx\pgfplots@listindex@set\pgfutil@empty
    \else
      \c@pgf@counta=\pgfplots@listindex
      \c@pgf@countb=\pgfplots@listindex@set
      \advance\c@pgf@countb by -\c@pgf@counta
      \globaldefs=1\relax
      \edef\setshift{%
        \noexpand\pgfkeys{
          /pgfplots/cycle list shift=\the\c@pgf@countb,
          /pgfplots/cycle list set=
        }
      }%
      \setshift%
      \globaldefs=0\relax
    \fi
    \pgfkeysgetvalue{/pgfplots/cycle list shift}\pgfplots@listindex@shift
    \ifx\pgfplots@listindex@shift\pgfutil@empty
    \else
      \c@pgf@counta=\pgfplots@listindex\relax
      \advance\c@pgf@counta by\pgfplots@listindex@shift\relax
      \ifnum\c@pgf@counta<0
        \c@pgf@counta=-\c@pgf@counta
      \fi
      \edef\pgfplots@listindex{\the\c@pgf@counta}%
    \fi
    \ifnum\pgfplots@cycle@dim>0
      \c@pgf@counta=\pgfplots@cycle@dim\relax
      \c@pgf@countb=\pgfplots@listindex\relax
      \advance\c@pgf@counta by-1
      \pgfplotsloop{%
        \ifnum\c@pgf@counta<0
          \pgfplotsloopcontinuefalse
        \else
          \pgfplotsloopcontinuetrue
        \fi
      }{%
        \pgfkeysgetvalue{/pgfplots/cycle multi list/@N\the\c@pgf@counta}\pgfplots@cycle@N
        \pgfplotsmathmodint{\c@pgf@countb}{\pgfplots@cycle@N}%
        \divide\c@pgf@countb by \pgfplots@cycle@N\relax
        \expandafter\pgfplots@getautoplotspec@
        \csname pgfp@cyclist@/pgfplots/cycle multi list/@list\the\c@pgf@counta @\endcsname
        {\pgfplots@cycle@N}%
        {\pgfmathresult}%
        \t@pgfplots@toka=\expandafter{#1,}%
        \t@pgfplots@tokb=\expandafter{\pgfplotsretval}%
        \edef#1{\the\t@pgfplots@toka\the\t@pgfplots@tokb}%
        \advance\c@pgf@counta by-1
      }%
    \else
      \pgfplotslistsize\autoplotspeclist\to\c@pgf@countd

      \pgfplots@getautoplotspec@{\autoplotspeclist}{\c@pgf@countd}{\pgfplots@listindex}%
      \let#1=\pgfplotsretval
    \fi
    \pgfmath@smuggleone#1%
    \endgroup
  }

\pgfplotsset{
  cycle list set/.initial=
}
\makeatother

\usepackage{lipsum}
\usepackage{amsfonts}
\usepackage{graphicx}
\usepackage{epstopdf}
\usepackage{amsmath,amssymb}
\ifpdf
  \DeclareGraphicsExtensions{.eps,.pdf,.png,.jpg}
\else
  \DeclareGraphicsExtensions{.eps}
\fi

\headers{Local Solvers for High-Order Patch Smoothers via p-Multigrid}{M. Wichrowski}

\title{Local Solvers for High-Order Patch Smoothers via p-Multigrid}

\author{MICHAŁ Wichrowski\thanks{Interdisciplinary Center for Scientifc Computing, Heidelberg University, Heidelberg,
    Germany}
}

\usepackage{amsopn}
\DeclareMathOperator{\diag}{diag}

\makeatletter
\newcommand*{\addFileDependency}[1]{
  \typeout{(#1)}
  \@addtofilelist{#1}
  \IfFileExists{#1}{}{\typeout{No file #1.}}
}
\makeatother

\begin{document}

\maketitle

\begin{abstract}
  I propose a vertex patch smoother where local problems are solved inexactly by a nested,
matrix-free \textit{p}-multigrid, creating a multigrid-within-multigrid framework. A single iteration of the
local solver can be evaluated with  $\mathcal{O}(p^{d+1})$ operations, and the approach is applicable to
non-separable problems on unstructured meshes. Numerical experiments demonstrate limited sensitivity
to geometric distortion and high-contrast coefficients. When used in a global geometric multigrid solver,
the method achieves robustness with respect to both polynomial degree $p$ and mesh refinement,
even on heavily distorted meshes.

\end{abstract}

\begin{keywords}
  finite element method, multigrid method, vertex-patch smoothing, data locality
\end{keywords}

\begin{AMS}
  65K05, 
  65Y10, 
  65Y20, 
  65N55, 
  65N30 
\end{AMS}

\section{Introduction}

Given their remarkable robustness with respect to the polynomial degree~\cite{pavarino1993additive} and excellent
data locality~\cite{wichrowski2025smoothers}, it is unsurprising that patch-based smoothers are gaining popularity in
modern high-performance finite element solvers. However, this power comes at a price: the core of a patch smoother is a
local solver applied to a small subdomain, or \emph{patch}, of the mesh. The efficiency of the overall method hinges on
the ability to solve these local problems rapidly. This has proven to be a persistent challenge, often leading to
complex, problem-specific solutions or a retreat to simpler, less effective smoothers.

In this paper, we embrace the recursive spirit of multigrid methods to address this challenge. Rather than treating the
local patch problems as a separate issue to be solved by external means, we apply the recursive logic of multigrid one
level deeper. We propose constructing the local solvers themselves using another multigrid solver --- this time,
employing p-coarsening. This \emph{multigrid-within-multigrid} approach yields a powerful, matrix-free local solver
that remains true to the iterative, computationally lean spirit of the global method. It avoids the need for direct
factorizations or other specialized techniques, offering a flexible and elegant solution to the local solver conundrum.

Schwarz smoothers can be viewed as a form of small-grain, overlapping domain
decomposition~\cite{pavarino1993additive}. The global problem is broken
down into a collection
of smaller, more manageable local problems defined on subdomains, or \emph{patches}, of the mesh. Next, the solutions
are pieced back together to form an improved global approximation. This approach has proven
highly effective, particularly for high-order discretizations and challenging problems like the Stokes
equations~\cite{hong2016robust,KanschatMao15, voronin2025monolithic}. Moreover, their applicability extends to non-matching grid methods, with
recent work demonstrating their effectiveness for the Shifted Boundary Method~\cite{wichrowski2025geometric,
  wichrowski2025matrix} and CutFEM~\cite{bergbauer2025high, cui2025multigrid}. A popular and simpler variant uses
non-overlapping,
cell-wise patches,
which can be seen as a form of block-Jacobi smoother~\cite{fischer2000overlapping, pazner2018approximate,
  remacle2016gpu}. Cell-wise smoothers have often been favored due to their simpler implementation and smaller data
structures for the local problems. However, this work challenges that trade-off: by introducing a matrix-free,
memory-optimal local solver, we make the more powerful vertex-patch smoothers just as lean and efficient, leveling the
playing field.

The efficiency of patch-based methods, especially for high-order discretizations, is tightly linked to matrix-free
computation. This reflects a broader shift toward matrix-free techniques that are essential for high-order finite
element methods~\cite{Kronbichler2012,Kronbichler2017a, kronbichler2019multigrid}. Traditional assembled sparse-matrix
approaches become prohibitively expensive
as the polynomial degree $p$ increases: the number of non-zero entries per row in the matrix grows like
$\mathcal{O}(p^d)$, leading to massive memory requirements that make the computation heavily memory-bandwidth-bound.
Matrix-free methods circumvent this bottleneck by never forming the global matrix, instead evaluating the operator's
action on-the-fly. For methods on tensor-product cells, this is particularly effective: sum-factorization techniques
reduce the operator application cost to $\mathcal{O}(p^{d+1})$, which is not only memory-optimal but also
asymptotically cheaper than a conventional sparse matrix-vector product~\cite{Kronbichler2017a}. Since high-order
methods are often
memory-bound, performance hinges on efficient data movement, making the excellent data locality of patch-based
approaches key to their success on modern hardware~\cite{munch2023cache, kronbichler2023enhancing, kronbichler2022cg}.
This has led to
highly performant matrix-free implementations on
both CPUs~\cite{wichrowski2025smoothers} and GPUs~\cite{cui2025implementation}. Our approach is designed to operate
entirely within this matrix-free paradigm, preserving its benefits at both the global and local levels.

While matrix-free techniques provide an efficient computational framework, they do not eliminate the challenge of
solving
the local patch problems. Several sophisticated methods have been proposed to tackle this, often diverging from a
purely iterative, matrix-free ideal. For instance, the work in~\cite{brubeck2021scalable} introduces a scalable direct
solver for separable problems on tensor-product cells. It constructs a special basis to diagonalize parts of the local
stiffness matrix, resulting in a sparser structure that makes Cholesky factorization feasible even for high polynomial
degrees. Alternatively, fast diagonalization techniques~\cite{WitteArndtKanschat21} offer an efficient way to invert
an operator on Cartesian reference patches. Despite their power, these methods rely on the
separability of the problem, restricting their use to structured grids with constant coefficients, and their
performance is known to deteriorate on deformed grids.

The complexity of vertex-patch solvers has led some to favor simpler alternatives~\cite{margenberg2025hp,
  anselmann2024energy}. A cell-wise
smoother, a form of block-Jacobi, restricts local problems to single cells, sometimes with a minimal overlap of degrees
of freedom from neighboring cells~\cite{fischer2000overlapping, pazner2018approximate, remacle2016gpu}. This simplifies
parallelization, as the local solves are entirely independent, but at the cost of reduced smoothing effectiveness
compared to larger, overlapping patches.

The cell-wise decomposition is naturally suited to iterative methods for the local solves, which aligns with the
matrix-free philosophy by avoiding matrix assembly altogether. For example, iterative solvers have been used for the
block-diagonal systems that arise in cell-based smoothers~\cite{bastian2019matrix}. While this shares some of the
spirit of our work, it is typically applied in a simpler, non-overlapping context and employs basic
diagonal preconditioners. Interestingly, for high-order discontinuous Galerkin methods, it has been shown that
transferring the problem to a continuous finite element space on the finest level can lead to robust and efficient
multigrid performance~\cite{fehn2020hybrid}.

Our work bridges this gap by introducing a p-multigrid method as a dedicated, high-performance iterative solver for the
local patch problems. This approach inherits the advantages of being matrix-free and algorithmically consistent with
the outer geometric multigrid method. It directly challenges the notion that robustness on challenging geometries
requires exact solvers. In particular, in~\cite{brubeck2021scalable} it was questioned whether p-robustness could be
maintained on highly distorted meshes even with exact patch solvers. Our results, summarized in
Tables~\ref{tab:gmres_iterations_cartesian_reinforced} and~\ref{tab:hp-comparison}, demonstrate the contrary: even with an approximate iterative solver on the patch, we achieve near-p-robustness on heavily distorted meshes, while the iteration counts decrease with mesh refinement.

The proposed method not only maintains optimal memory complexity of $\mathcal{O}(p^d)$ and requires minimal
precomputation but also proves resilient to large, discontinuous jumps in material coefficients. On top of that,
our patch-based approach, by its nature, provides excellent data locality~\cite{wichrowski2025smoothers}.

While we validate the method on the Laplacian, our idea is to provide a flexible framework for the local
patch-solve problem. The p-multigrid solver presented here is highly effective, but it could be further improved or
even replaced by other advanced methods, such as the direct solver from~\cite{brubeck2021scalable} on Cartesian-like
patches where it excels~\cite{WitteArndtKanschat21}. Furthermore, the concept of a recursive multigrid-based local
solver is general. It can be
extended to more complex systems, such as the Stokes equations, for which p-multigrid solvers have already been
developed~\cite{braess1997efficient,wichrowski2022matrix, jodlbauer2024matrix}. We thus present a practical and
efficient smoother that extends the applicability of patch-based multigrid to a wider range of challenging problems.
The implementation is done using the open-source finite element library \texttt{deal.II}~\cite{dealII97,dealii2019design}.

The remainder of this paper is organized to guide the reader from the foundational concepts to the final numerical
results. In Section~\ref{sec:multigrid} we briefly describe the
multigrid method and its key ingredient, the smoothing step. The core idea of this paper, the p-multigrid local solver,
is detailed in Section~\ref{sec:local_solve}. Finally, Section~\ref{sec:application_in_mg} presents numerical results
that validate the proposed method, with key performance indicators summarized in
Table~\ref{tab:gmres_iterations_cartesian_reinforced} for a quick overview of its effectiveness.

\section{Geometric finite element multigrid methods}
\label{sec:multigrid}
Let us first briefly describe the mathematical setup of the multigrid method, first
assume a hierarchy of meshes
\begin{gather}
  \mesh_1 \sqsubset \mesh_2 \sqsubset \dots \sqsubset \mesh_L,
\end{gather}
subdividing a domain in $\R^d$, where the symbol
``$\sqsubset$'' indicates nestedness, that is, every cell of mesh
$\mesh_{\ell+1}$ is obtained from a cell of mesh $\mesh_{\ell}$ by refinement. We
note that topological nestedness is sufficient from the algorithmic
point of view, such that domains with curved boundaries can be
covered approximately.

With each mesh $\mesh_\ell$ we associate a finite element space $V_\ell$
spanned by locally supported shape functions with their associated degrees of freedom.
We identify $V_\ell$ with $\mathbb R^{\operatorname{dim}V_\ell}$ and do not distinguish between a finite element
function $u_\ell$ and its coefficient vector by notation.
Between these spaces, we introduce transfer operators
\begin{gather}
  \begin{array}{rlcl}
    \resth\ell\colon & V_{\ell+1} & \to & V_\ell,     \\
    \prolh\ell\colon & V_{\ell}   & \to & V_{\ell+1}.
  \end{array}
\end{gather}
As usual, $\prolh\ell$ is chosen as the embedding operator and $\resth\ell$ as its
$\ell_2$-adjoint. The subscript $h$ indicates h-multigrid transfer between meshes of different resolution.

For our numerical experiments, we consider the variable-coefficient Poisson equation with bilinear form
\begin{gather}
  \label{eq:bilinear-form}
  a(u,v) = \int_\Omega \mu \nabla u \cdot \nabla v \, d\vec x,
\end{gather}
where $\mu$ is a spatially varying coefficient.
The discrete problem on level $\ell$ reads
\begin{gather}
  \label{eq:matrix}
  A_\ell u_\ell = b_\ell,
\end{gather}
where $A_\ell$ and $b_\ell$ are obtained by applying $a(\cdot,\cdot)$ and $b(\cdot)$ to the basis functions of
$V_\ell$.

Multigrid methods are efficient solution methods for the discrete linear system \eqref{eq:matrix}, at least for PDE
with a regularity gain. When implementing the finite element V-cycle, see for instance~\cite{Bramble93,Hackbusch85}, a
so-called smoother is needed on each level in addition to the operators $A_\ell$, $\resth\ell$, and $\prolh\ell$. It is
usually described by its error propagation operator $\mathcal{S_{\ell}}$.
The simplest example is Richardson's method, where
\begin{gather*}
  \mathcal{S_{\ell,\text{Rich}}} = I_\ell - \omega A_\ell.
\end{gather*}
We introduce more powerful smoothers in the next subsection.

\subsection{Vertex patch smoothers}
\label{sec:subspace-correction}

The vertex patch smoother is an overlapping subspace correction
method~\cite{ArnoldFalkWinther00,JanssenKanschat11,KanschatMao15,WitteArndtKanschat21,brubeck2021scalable,wichrowski2025smoothers}.
The domain
is decomposed into a collection of overlapping subdomains, or \emph{patches} $\Omega_j$, each formed by the cells
surrounding an interior vertex $x_j$. The smoother then iteratively improves the solution by solving local problems on
these patches. We employ a multiplicative (Gauss-Seidel-alike) variant, where patches are processed sequentially and
the
global solution is updated after each local solve. To avoid recomputing the full global residual for every patch,
we use a technique similar to the one proposed in~\cite{wichrowski2025smoothers}.

The update for a single patch $j$ consists of four steps, which directly correspond to the practical implementation of
the local correction operator $P_j v = \Pi_j ^T \tilde A_j^{-1} \Pi_j (b-Au)$. Here, $\Pi_j$ extracts local data from
global vectors, and $\tilde A_j^{-1}$ is an approximate inverse of the local operator $A_j$. The steps are as follows:
\begin{enumerate}
  \item \textbf{Gather.} Collect the current solution values $u$ and right-hand side values $b$ corresponding to the
        degrees of freedom on the patch $\Omega_j$ and its boundary. We denote the spaces of interior and boundary
        degrees of
        freedom by $V_j$ and $\overline{V}_j$, respectively, and the gathering operators by $\Pi_j$ and
        $\overline{\Pi}_j$.
        \begin{gather}
          \label{eq:gather}
          \overline u_j = \overline\Pi_j u \in \overline V_j, \qquad b_j = \Pi_j b \in V_j
        \end{gather}
  \item \textbf{Evaluate.} Compute the local residual $r_j$ for the interior degrees of freedom using the gathered data
        and the local operator $\overline{A}_j$.
        \begin{gather}
          \label{eq:main:4}
          r_j = b_j - \Pi_j\overline A_j \overline u_j \in V_j.
        \end{gather}
  \item \textbf{Local solve (LS).} Find a correction $d_j$ by approximately solving the local system on the patch
        interior.
        \begin{gather}
          \label{eq:main:5}
          d_j \approx \tilde A_j^{-1} r_j,
        \end{gather}
        where $\tilde A_j^{-1}$ is an inexpensive approximate solver, which is the main topic of this paper.
  \item \textbf{Scatter.} Add the local correction $d_j$ back to the global solution vector $u$.
        \begin{gather}
          \label{eq:scatter}
          u \pluseq \Pi_j^T d_j.
        \end{gather}
\end{enumerate}

\begin{figure}[tp]
  \centering
  \def\svgwidth{\columnwidth}
  \begin{tikzpicture}[scale=0.8, transform shape, font=\normalsize,
      cell/.style={draw=black, thick, minimum size=1.2cm},
      interior/.style={fill=blue!20},
      global/.style={fill=red!15},
      panel-label/.style={font=\bfseries\small, align=center},
      arrow/.style={->, thick, shorten >=2mm, shorten <=2mm},
      dof/.style={circle, fill=black, inner sep=0.8pt}
    ]

    \newcommand{\drawPatchInteriorDoFs}[2]{
      \foreach \i in {1,...,5} {
          \foreach \j in {1,...,5} {
              \pgfmathsetmacro{\xi}{\i/6}
              \pgfmathsetmacro{\eta}{\j/6}

              \coordinate (NW) at (#1.north west);
              \coordinate (SE) at (#2.south east);

              \coordinate (pt_x) at ($(NW)!\xi!(SE)$);   
              \coordinate (pt_y) at ($(NW)!\eta!(SE)$);  
              \coordinate (pt) at (pt_x |- pt_y);
              \node[dof, minimum size=4pt, inner sep=0pt, fill=blue] at (pt) {};
            }
        }
    }

    \newcommand{\drawCubicDoFs}[1]{
      \foreach \xi in {0, 0.33, 0.66, 1} {
          \foreach \eta in {0, 0.33, 0.66, 1} {
              \coordinate (bottom) at ($(#1.south west)!\xi!(#1.south east)$);
              \coordinate (top) at ($(#1.north west)!\xi!(#1.north east)$);
              \node[dof, minimum size=3pt, inner sep=0pt,  fill=red] at ($(bottom)!\eta!(top)$) {};
            }
        }
    }

    \newcommand{\drawBoundaryDoFs}[2][dof]{
      \foreach \xi in {0, 0.33, 0.66, 1} {
          \coordinate (pt) at ($(#2.south west)!\xi!(#2.south east)$);
          \node[#1] at (pt) {};
        }
      \foreach \xi in {0, 0.33, 0.66, 1} {
          \coordinate (pt) at ($(#2.north west)!\xi!(#2.north east)$);
          \node[#1] at (pt) {};
        }
      \foreach \eta in {0.33, 0.66} {
          \coordinate (pt) at ($(#2.south west)!\eta!(#2.north west)$);
          \node[#1] at (pt) {};
        }
      \foreach \eta in {0.33, 0.66} {
          \coordinate (pt) at ($(#2.south east)!\eta!(#2.north east)$);
          \node[#1] at (pt) {};
        }
    }
    \newcommand{\drawInteriorFrame}[2]{%
      \draw[rounded corners=3pt, thick, draw=blue!60!black, fill=gray!80, fill opacity=0.1, densely dashed]
      ([shift={(4pt,-4pt)}]#1.north west)
      rectangle
      ([shift={(-4pt,4pt)}]#2.south east);%
    }

    \node[cell] (p1c1) at (0,0) {};
    \node[cell] (p1c2) at (1.2,0) {};
    \node[cell] (p1c3) at (0,-1.2) {};
    \node[cell] (p1c4) at (1.2,-1.2) {};
    \drawInteriorFrame{p1c1}{p1c4}

    \drawCubicDoFs{p1c1}
    \drawCubicDoFs{p1c2}
    \drawCubicDoFs{p1c3}
    \drawCubicDoFs{p1c4}
    \drawPatchInteriorDoFs{p1c1}{p1c4}
    \node[panel-label, below=0.3cm of p1c3.south west, anchor=north west, ] {Solution\\ Right-hand side\\(all DoFs)};

    \node[cell] (p2c1) at (4.5,0) {}; 
    \node[cell] (p2c2) at (5.7,0) {}; 
    \node[cell] (p2c3) at (4.5,-1.2) {}; 
    \node[cell] (p2c4) at (5.7,-1.2) {}; 
    \drawInteriorFrame{p2c1}{p2c4}
    \drawPatchInteriorDoFs{p2c1}{p2c4}
    \node[panel-label, below=0.3cm of p2c3.south west, anchor=north west, ] { Residual \\ (interior DoFs)};

    \node[cell] (p3c1) at (9,0) {}; 
    \node[cell] (p3c2) at (10.2,0) {}; 
    \node[cell] (p3c3) at (9,-1.2) {}; 
    \node[cell] (p3c4) at (10.2,-1.2) {}; 
    \drawInteriorFrame{p3c1}{p3c4}
    \drawPatchInteriorDoFs{p3c1}{p3c4}
    \node[panel-label, below=0.3cm of p3c3.south west, anchor=north west,] {Correction\\(interior only)};

    \draw[arrow] ($(p1c3.west) - (2cm, 0)$) --  node[above, midway, ] {$\overline\Pi_j u, \;\;  \Pi_j b $} node[below, midway, ] {Gather}  (p1c3.west);
    \draw[arrow] (p1c4.east) -- node[above, midway, ] {$b_j - \Pi_j\overline A_j \overline u_j$} node[below, midway, ] {Evaluate} (p2c3.west);
    \draw[arrow] (p2c4.east) -- node[above, midway, ]{$\tilde A_j^{-1} r_j$}
    node[below, midway, ] {Local Solve}   (p3c3.west);
    \draw[arrow] (p3c4.east) -- node[above, midway, ]{$\Pi_j^T d_j$} node[below, midway, ] {Scatter} ($(p3c4.east) + (2cm, 0)$);

  \end{tikzpicture}

  \caption{Workflow in smoother application on a single patch}
  \label{fig:patch-smoother-steps}
\end{figure}
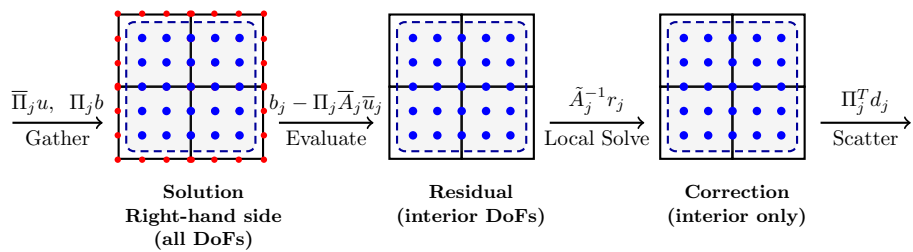

These steps are visualized in Figure~\ref{fig:patch-smoother-steps} and summarized in
Algorithm~\ref{alg:Loop-sequential} together with a smoother sweep. This implementation follows the cell-oriented loop
structure proposed
in~\cite{wichrowski2025smoothers}, which is well-suited for matrix-free operator evaluation. Building on this
framework, we can employ iterative solvers for the local problem, avoiding dense direct factorizations.

\begin{algorithm}[tp]

  \begin{algorithmic}

    \State \Function{\LUpNew}{j,u}
    \State    $ u_j \gets \overline\Pi_j u$ \Comment*{Gather}
    \State    $ r_{j} \gets \Pi_j b - \Pi_j\overline A_j u_j$ \Comment*{Evaluate}
    \State    $ d_j \gets \tilde{A}_j ^{-1} r_j$ \Comment*{Solve}
    \State \Return{ $\Pi_j^T d_j$ }\Comment*{Scatter}
    \EndFunction

    \State \For{$j=1,\ldots,N_\text{patches}$}
    \State $u \pluseq$ \Call{\LUpNew}{$j$,$u$}

    \EndFor

  \end{algorithmic}
  \caption{Application of a patch smoother.}
  \label{alg:Loop-sequential}
\end{algorithm}

\section{Local solvers via p-multigrid}
\label{sec:local_solve}
The development of multigrid methods can be motivated by the question of how to accelerate simple iterative solvers
like Jacobi or Gauss-Seidel. While these smoothers are effective at reducing high-frequency error components, they are
slow to converge for low-frequency components. The introduction of another level of coarse-grid correction forms a
two-level method, which addresses this by efficiently handling these low-frequency errors on a coarser grid. This
naturally leads to the next question: how should the coarse-grid problem be solved? The most effective approach is to
again introduce another level, recursively leading to the full multigrid algorithm.

When considering the local solve on each patch (Equation~\eqref{eq:main:5}), we face a similar challenge: we need an
efficient solver for a smaller, but still potentially non-trivial, linear system. Following the same logic that led to
the development of multigrid for the global problem, a natural and consistent choice for the local solver is to employ
additional levels. We therefore treat each patch solve as a small multigrid problem in its own right.

For geometric (h-)multigrid, this recursion is restricted by the limited number of refinement levels available inside
a single patch. Since a patch consists of only $2^d$ cells, one coarsening step already leads to a single cell, which
is problematic for several reasons. For instance, for linear elements, no degrees of freedom would be left after
imposing boundary conditions on the patch boundary. Furthermore, for some problems like the Stokes equations, smoothing
over multiple cells is necessary for convergence~\cite{hong2016robust}. To overcome these limitations, we employ
p-multigrid (coarsening in polynomial degree) as the inner solver on each patch: the geometry and mesh of the patch
remain fixed while we recursively reduce the polynomial degree of the local approximation. This choice preserves the
advantages of the
matrix-free, sum-factorized operator evaluation provided by the local framework, yields inexpensive coarse solves and
keeps iteration counts effectively bounded with respect to the polynomial degree.

\subsection{p-multigrid components}
\label{sec:p_mg_components}

To construct the p-multigrid solver for the local patch problem \eqref{eq:main:5}, we define a hierarchy of finite
element spaces on the patch, indexed by polynomial degree $p$:
\begin{gather}
  V_{j,1} \subset V_{j,2} \subset \dots \subset V_{j,p_{\text{max}}},
\end{gather}
where $V_{j,p}$ is the space of polynomials of degree up to $p$ on the patch $\Omega_j$. As for any multigrid
algorithm, the components of our p-multigrid cycle are the transfer operators, the coarse-level operators, and the
smoothers.

\paragraph{Transfer operators}
The transfer operators between different polynomial degrees are defined by standard $L^2$-projections. The prolongation
operator $\prolp{p-1,p}\colon V_{j,p-1} \to V_{j,p}$ is the natural embedding, while the restriction operator
$\restp{p,p-1}\colon V_{j,p} \to V_{j,p-1}$ is its adjoint.

The transfer operators inherit a separable (Kronecker) structure $\mathcal R^{\otimes d}$ and can therefore be applied
by $d$ successive one-dimensional tensor contractions. Using sum-factorization this leads to an arithmetic cost scaling
like $\mathcal O(p^{d+1})$ (more precisely $\mathcal O(d\,p^{d+1})$).

\paragraph{Level operators}
Formally, the operators on the p-multigrid levels are defined by the Galerkin projection
\begin{gather}
  A_{j,p-1} = \restp{p,p-1}\, A_{j,p}\, \prolp{p-1,p}.
\end{gather}
This relation defines the coarse-level operator as the Galerkin projection of the fine-level operator. Since the
space $V_{j,p-1}$ is identified as a subspace of $V_{j,p}$ (and its basis is chosen accordingly), this projection
corresponds to selecting the sub-block / the representation of $A_{j,p}$ associated with the lower-degree basis.
In our implementation we therefore simply evaluate the operator on the coarse p-level using the basis functions of the
lower-degree finite element space.

The application of the level operators using matrix-free sum factorization incurs an arithmetic cost scaling like
$\mathcal{O}(p^{d+1})$ per cell, since there are $\mathcal{O}(p^d)$ quadrature points and each tensorized evaluation
requires $\mathcal{O}(p)$ work along one coordinate. Memory traffic, however, scales like $\mathcal{O}(p^d)$ for local
vectors and coefficients, so for large polynomial degrees the computation becomes increasingly compute-bound.

\paragraph{Coarse solver}
On the coarsest p-level ($p=1$) the patch interior contains only a single unconstrained node. Hence  for a scalar
problem the local operator $A_{j,1}$ is a $1\times1$ matrix and its inverse is just the reciprocal of that scalar:
\begin{gather*}
  A_{j,1}^{-1} = \frac{1}{A_{j,1}}.
\end{gather*}
In practice we store this single number and apply it as the coarse-level solver. Note that this scalar equals the
diagonal entry of the operator assembled on the patch with polynomial degree one (or, equivalently, the corresponding
diagonal entry of the fine-level operator after Galerkin projection), so its extraction and storage are inexpensive.

\paragraph{Smoothers}
For the smoothing step within the p-multigrid cycle, we use a preconditioned Richardson iteration with   a simple preconditioner
$P_p$ and  a damping parameter $\omega$. A computationally inexpensive choice is a damped Jacobi smoother, for
which $P_p = D_{j,p}$ is the inverse diagonal of the local matrix $A_{j,p}$.

Another powerful approach for Cartesian grids is the fast diagonalization method and use it to construct a smoother.
The idea is to use the inverse of the operator on a reference Cartesian patch, $A_{b,p}$, scale it with a diagonal
matrix, and use it as a preconditioner inside the smoothing step. The exact form of the preconditioner is
\begin{gather*}
  P_p = D_{j,p} \diag(A_{b,p})  A_{b,p}^{-1}.
\end{gather*}
where $A_{b,p}$ is the operator on the reference Cartesian patch with polynomial degree $p$.
The diagonal scaling with $D_{j,p} \diag(A_{b,p})$ is  supposed to adjust for variations in geometry and coefficients.
This approach uses the efficient structure of the reference problem to
accelerate smoothing for more general cases, however is limited to structured patches.

The complete p-multigrid V-cycle algorithm, used to approximately solve the local patch problem $A_{j,p} d_p = r_p$, is
summarized in Algorithm~\ref{alg:p-v-cycle}.
\begin{algorithm}[htbp]
  \caption{The p-multigrid V-cycle for a local patch problem.}
  \label{alg:p-v-cycle}
  \begin{algorithmic}[1]
    \Function{V-cycle}{$p, d_p, r_p$}
    \If{$p = 1$}
    \State $d_1 \gets A_{j,1}^{-1} r_1$ \Comment{Coarse-level solve}
    \State \Return $d_1$
    \EndIf
    \State $d_p \gets d_p + P_p^{-1} (r_p - A_{j,p} d_p)$ \Comment{Pre-smoothing}
    \State $res_p \gets r_p - A_{j,p} d_p$ \Comment{Compute residual}
    \State $res_{p-1} \gets \restp{p,p-1} res_p$ \Comment{Restrict residual}
    \State $e_{p-1} \gets \Call{V-cycle}{p-1, 0, res_{p-1}}$ \Comment{Recursive call}
    \State $e_p \gets \prolp{p-1,p} e_{p-1}$ \Comment{Prolongate correction}
    \State $d_p \gets d_p + e_p$ \Comment{Apply correction}
    \State $d_p \gets d_p + P_p^{-1} (r_p - A_{j,p} d_p)$ \Comment{Post-smoothing}
    \State \Return $d_p$
    \EndFunction
  \end{algorithmic}
\end{algorithm}

Estimating the cost of the smoother is straightforward once one recognises that the dominant operation on each patch is
the matrix-free operator evaluation implemented by sum-factorisation. To compute the local residual one must evaluate
gradients at quadrature points and then integrate the resulting fluxes back to the coefficient space. For
tensor-product elements this requires $d$ one-dimensional sum-factorisations to form the directional derivatives (one
per
coordinate) and another $d$ one-dimensional sum-factorisations to integrate those contributions, i.e. roughly $2d$
one-dimensional sweeps. Consequently a single matrix-free operator application costs on the
order of $O(p^{\,d+1})$ floating point operations.

As a results, fast diagonalization and the p-multigrid inner solver share the same asymptotic work for a single,
matrix-free
operator application: both approaches are dominated by tensorized (sum-factorized) contractions and therefore incur
roughly $O(d\,p^{d+1})$ floating point operations on a patch of polynomial degree p in d dimensions. The crucial
difference is how many such applications are required.

\subsection{Numerical validation on a single patch}
\label{sec:local_solve_validation}
Before integrating the local p-multigrid solver into the global smoothing procedure, we first evaluate its performance.
For these tests, we run
the CG solver using a p-multigrid V-cycle as a preconditioner until the relative residual is reduced below a tolerance
of $10^{-8}$. In all experiments,
we choose the damping parameter $\omega = \nicefrac{1}{2}$. We emphasize that such a stringent tolerance is not
necessary when the p-multigrid is used as a patch smoother in the global algorithm---a much looser tolerance typically
suffices and is computationally cheaper. We nevertheless enforce $10^{-8}$ here to increase the number of CG
iterations,
which makes convergence behavior easier to inspect and provides more detailed diagnostics for validation.
These tests are not intended to validate the local solver as a standalone production solver; rather, they isolate its
behaviour from other multigrid components and assess only the contributions proposed in this paper.

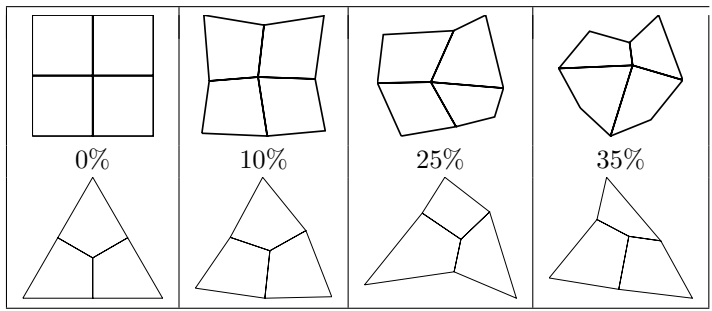
\begin{figure}[htbp]
  \centering
\begin{tabular}{|c|c|c|c|}
  \hline
      &      &             \\[-0.9em]
  \begin{tikzpicture}[scale=0.8, every node/.style={scale=0.8}]
    \begin{axis}[
        axis equal image,
        hide axis,
        width=0.32\columnwidth,
        enlargelimits=false
      ]
      \addplot[black,thick,mark=none] table {figures/patch_def_data/patch00_2D_2.gp};
    \end{axis}
  \end{tikzpicture}
      &
  \begin{tikzpicture}[scale=0.8, every node/.style={scale=0.8}]
    \begin{axis}[
        axis equal image,
        hide axis,
        width=0.32\columnwidth,
        enlargelimits=false
      ]
      \addplot[black,thick,mark=none] table {figures/patch_def_data/patch01_2D.gp};
    \end{axis}
  \end{tikzpicture}
      &
  \begin{tikzpicture}[scale=0.8, every node/.style={scale=0.8}]
    \begin{axis}[
        axis equal image,
        hide axis,
        width=0.32\columnwidth,
        enlargelimits=false
      ]
      \addplot[black,thick,mark=none] table {figures/patch_def_data/patch025_2D_1.gp};
    \end{axis}
  \end{tikzpicture}
      &
  \begin{tikzpicture}[scale=0.8, every node/.style={scale=0.8}]
    \begin{axis}[
        axis equal image,
        hide axis,
        width=0.32\columnwidth,
        enlargelimits=false
      ]
      \addplot[black,thick,mark=none] table {figures/patch_def_data/patch035_2D.gp};
    \end{axis}
  \end{tikzpicture}
  \\
  0\% & 10\% & 25\% & 35\%
  \\
  \begin{tikzpicture}[scale=0.8, every node/.style={scale=0.8}]
    \begin{axis}[
        axis equal image,
        hide axis,
        width=0.32\columnwidth,
        enlargelimits=false
      ]
      \addplot[black,mark=none] table {figures/patch_def_data/simplex00_2D_2.gp};
    \end{axis}
  \end{tikzpicture}
      &
  \begin{tikzpicture}[scale=0.8, every node/.style={scale=0.8}]
    \begin{axis}[
        axis equal image,
        hide axis,
        width=0.32\columnwidth,
        enlargelimits=false
      ]
      \addplot[black,mark=none] table {figures/patch_def_data/simplex01_2D.gp};
    \end{axis}
  \end{tikzpicture}
      &
  \begin{tikzpicture}[scale=0.8, every node/.style={scale=0.8}]
    \begin{axis}[
        axis equal image,
        hide axis,
        width=0.32\columnwidth,
        enlargelimits=false
      ]
      \addplot[black,mark=none] table {figures/patch_def_data/simplex025_2D.gp};
    \end{axis}
  \end{tikzpicture}
      &
  \begin{tikzpicture}[scale=0.8, every node/.style={scale=0.8}]
    \begin{axis}[
        axis equal image,
        hide axis,
        width=0.32\columnwidth,
        enlargelimits=false
      ]
      \addplot[black,mark=none] table {figures/patch_def_data/simplex035_2D.gp};
    \end{axis}
  \end{tikzpicture}
  \\
  \hline
\end{tabular}
\caption{Deformed structured (top row) and unstructured (bottom row) patches in 2D. The distortion parameter $\delta$
  is $0\%$, $10\%$, $25\%$, and $35\%$ (left to right).}
  \label{fig:deformed_patches}
\end{figure}
We validate the p-multigrid patch solver on a single vertex patch consisting of $2^d$ cells in spatial dimensions $d=2$
and $d=3$. In addition to Cartesian patches, we also test simplicial (unstructured) patches consisting of 3
quadrilaterals in 2D and 4 hexahedra in 3D; example patches are shown in Figure~\ref{fig:deformed_patches}. The
numerical study starts from a simple Poisson problem  with a constant coefficient on an undeformed
patch as a baseline, and then progressively increases difficulty to probe robustness and limitations.

First, we introduce geometric complexity by randomly displacing each vertex.
Each interior vertex is displaced by a vector of length $\delta h$, with its direction sampled uniformly at random
(i.e., uniformly on the unit sphere $S^{d-1}$).
The effect of this perturbation is illustrated in Figure~\ref{fig:deformed_patches}.

For larger distortions some generated patches may contain
non‑shape‑regular (degenerate) elements, e.g., elements with concave or inverted corners. For a 2D Cartesian patch,
degenerate elements occur when the distortion exceeds $\frac{1}{2}\frac{\sqrt{2}}{2} \approx 0.3535$, while in 3D the
threshold is $\frac{1}{2}\frac{\sqrt{3}}{3} \approx 0.2887$. Accordingly, we test distortions up to 35\% in 2D and up
to 30\% in 3D. For simplicial patches, degeneracy occurs at lower distortion
levels because the initial angles already
exceed $90^\circ$; such degenerate patches are identified and excluded from our experiments. Note that in
Figure~\ref{fig:deformed_patches}, some cells in the structured patch are close to triangles, with some angles
approaching $180^\circ$.

Second, we investigate the influence of heterogeneous
material properties, focusing on discontinuous jumps where one cell carries a large contrast relative to the others.
Finally, we combine geometric deformation and strong coefficient variation to evaluate the method under the most
challenging conditions. For each scenario, we run experiments across a range of polynomial degrees, use a single
p-multigrid V-cycle as a preconditioner within a CG iteration, and report iteration counts to quantify performance with
respect to $p$, mesh distortion, and coefficient contrast. Since the results in 2D and 3D are qualitatively very
similar, and the computations using sparse matrices are slow, we focus most experiments on the 2D case and present 3D
results only for selected baseline scenarios.

We perform 20 independent random realizations when testing geometric distortion and report the average number of
iterations over all runs. Each local solve is allowed a maximum of 100 iterations; if a run does not reach the
convergence tolerance within this limit, it is declared non-convergent. Such runs are assigned a large sentinel
iteration count, so that the computed average is driven well beyond the typical plot limits and the corresponding point
is visibly off the axis, making non-convergence immediately apparent in the figures.

\begin{figure}[htbp]
  \centering
  \usepgfplotslibrary{fillbetween}
  \begin{tikzpicture}
    \begin{axis}[
        longplot,
        xlabel={polynomial degree $p$},
        ylabel={CG iterations},
        xmin=3, xmax=16,
        grid=major,
        legend pos=north west,
        cycle list name=paper markers,
        every axis/.append style={font=\small},
        thick
      ]
      \addplot+[ thick] table[
          col sep=comma,
          x=degree,
          y=avg,
        ] {results/Jac_cart_2D.csv};
      \addlegendentry{$2D$}

      \addplot+[ thick] table[
          col sep=comma,
          x=degree,
          y=avg,
        ] {results/Jac_cart_3D.csv};
      \addlegendentry{$3D$}

      \pgfplotsset {cycle list set=0}
      \addplot+[thick, dashed] table[
          col sep=comma,
          x=degree,
          y=avg,
        ] {results/simplex_const_p_2D.csv};

      \addplot+[thick, dashed] table[
          col sep=comma,
          x=degree,
          y=avg,
        ] {results/simplex_const_p_3D.csv};

    \end{axis}
  \end{tikzpicture}
  \caption{Number of CG iterations versus polynomial degree $p$ for solving the local patch problem on a Cartesian
    patch in 2D and 3D. The CG method is preconditioned with a single p-multigrid V-cycle; the dashed line denotes
    results for a non-structured (simplicial) patch.}
  \label{fig:iter_vs_p_Jacobi_cartesian}
\end{figure}
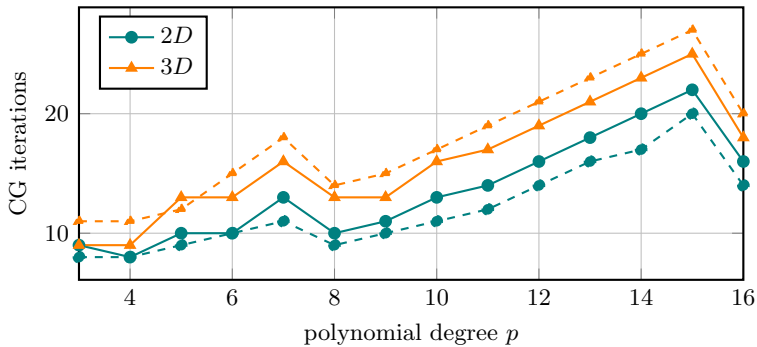

\paragraph{Baseline Performance}
We begin by evaluating the baseline performance of the p-multigrid preconditioner with a damped Jacobi smoother on
undeformed patches with a constant coefficient. Figure~\ref{fig:iter_vs_p_Jacobi_cartesian} reports the number of CG
iterations required for convergence as a function of the polynomial degree $p$. For undeformed Cartesian patches
(solid lines) using the Jacobi-based preconditioner, iteration counts increase moderately with $p$, from about 10 to 23
as $p$ grows from 3 to 15. The dashed line, showing results for a non-structured (simplex) patch in 2D, exhibits a
similar growth with $p$. We emphasize that these baseline results refer specifically to the Jacobi smoother: the
Cartesian-reinforced smoother yields exact local solves for the constant-coefficient, structured patches
considered here and therefore is not included in this particular comparison.

We note an apparent drop in iterations at $p=7$ and $p=15$. This is an artifact of the p-multigrid level structure: we
build p-levels using the nested sequence 1, 3, 7, 15, etc. (i.e. $p_{k+1} = 2p_k + 1$), so degrees that
coincide with these coarsening nodes benefit from a different multigrid hierarchy. For example, when $p=8$, the actual
V-cycle levels are \{1, 3, 7, 8\}, which effectively places the top two degrees very close and changes how smoothing
and coarse correction are applied. Such level-structure effects cause the non-monotonic behavior visible in the plot.

\begin{figure}[htbp]
  \begin{tikzpicture}
    \begin{scope}
      \begin{axis}[
          paperplot,
          xlabel={distortion},
          ylabel={CG iterations},
          xmin=0.0, xmax=0.35,
          xtick={0.05,0.1,0.15,0.2,0.25,0.3,0.35},
          ymin=0, ymax=25,
          legend pos=north west,
          axis lines=box,
          grid=major,
          clip=true,
        ]
        \node[anchor=north east,inner sep=5pt,fill=white,fill opacity=0.9,text opacity=1,draw=black,rounded
          corners=2pt] at (rel axis cs:0.9,0.24) {2D};

        \addplot+[thick] table[
            col sep=comma,
            x=distortion,
            y=avg,
            restrict expr to domain={\thisrow{degree}}{2.9:3.1},
          ] {results/cartesian_const_2D.csv};
        \addlegendentry{$p=3$}

        \addplot+[ thick] table[
            col sep=comma,
            x=distortion,
            y=avg,
            restrict expr to domain={\thisrow{degree}}{6.9:7.1},
          ] {results/cartesian_const_2D.csv};
        \addlegendentry{$p=7$}

        \addplot+[thick] table[
            col sep=comma,
            x=distortion,
            y=avg,
            restrict expr to domain={\thisrow{degree}}{14.9:15.1},
          ] {results/cartesian_const_2D.csv};
        \addlegendentry{$p=15$}

      \end{axis}
    \end{scope}

    \begin{scope}[xshift=0.52\columnwidth]
      \begin{axis}[
          paperplot,
          xlabel={distortion},
          ylabel={CG iterations},
          xmin=0.0, xmax=0.3,
          xtick={0.05,0.1,0.15,0.2,0.25,0.3,0.35},
          ymin=0, ymax=25,
          legend pos=north west,
        ]
        \node[anchor=north east,inner sep=5pt,fill=white,fill opacity=0.9,text opacity=1,draw=black,rounded
          corners=2pt] at (rel axis cs:0.9,0.24) {3D};

        \addplot+[ thick] table[
            col sep=comma,
            x=distortion,
            y=avg,
            restrict expr to domain={\thisrow{degree}}{2.9:3.1},
          ] {results/cartesian_const_3D.csv};
        \addlegendentry{$p=3$}

        \addplot+[ thick] table[
            col sep=comma,
            x=distortion,
            y=avg,
            restrict expr to domain={\thisrow{degree}}{6.9:7.1},
          ] {results/cartesian_const_3D.csv};
        \addlegendentry{$p=7$}

        \addplot+[ thick] table[
            col sep=comma,
            x=distortion,
            y=avg,
            restrict expr to domain={\thisrow{degree}}{14.9:15.1},
          ] {results/cartesian_const_3D.csv};
        \addlegendentry{$p=15$}

      \end{axis}
    \end{scope}
  \end{tikzpicture}
  \caption{Number of CG iterations vs.\ mesh distortion for the Cartesian-reinforced smoother.
    The CG solver is preconditioned with a single p-multigrid V-cycle applied as the Cartesian-reinforced patch
    solver. Each point reports the average number of CG iterations required to reduce the relative residual to
    $10^{-8}$. The left and right subplots show results for 2D and 3D vertex patches, respectively; the horizontal axis
    denotes the mesh distortion parameter used to deform
    the patches.}
  \label{fig:iter_vs_distortion}
\end{figure}
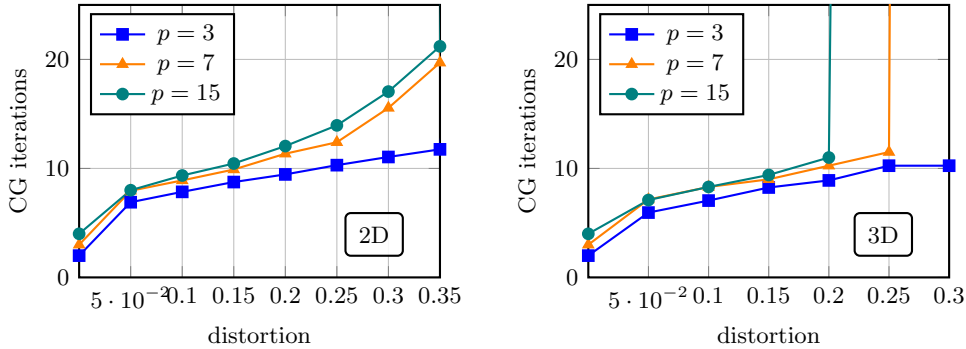

\paragraph{Resilience to Geometric Distortion}
Next, we examine the solver's performance under geometric perturbations. We apply a random, nonlinear distortion to the
patch geometry controlled by a \emph{distortion} parameter, with zero corresponding to an undeformed patch.
Figure~\ref{fig:iter_vs_distortion} presents the results for the Cartesian-reinforced smoother. This smoother is
designed to be effective for structured or near-structured grids by using an inverted operator on an ideal Cartesian
reference patch. The left and right panels show results for 2D and 3D, respectively. For small distortions, the
performance is exceptional, with very few CG iterations required for all tested polynomial degrees. As the distortion
increases, the iteration count grows, which is expected since the actual patch geometry deviates more from the
reference patch assumed by the smoother. In 2D, this growth is moderate, indicating that the smoother remains effective
even for reasonably deformed patches. However, in 3D, the method is more sensitive to distortion, especially for higher
polynomial degrees. For instance, with $p=15$, convergence issues arise for distortions around 20\%, and for $p=7$, the
limit is reached at about 25\% distortion. This indicates that high-order basis functions are more affected by
geometric perturbations, particularly in three dimensions.

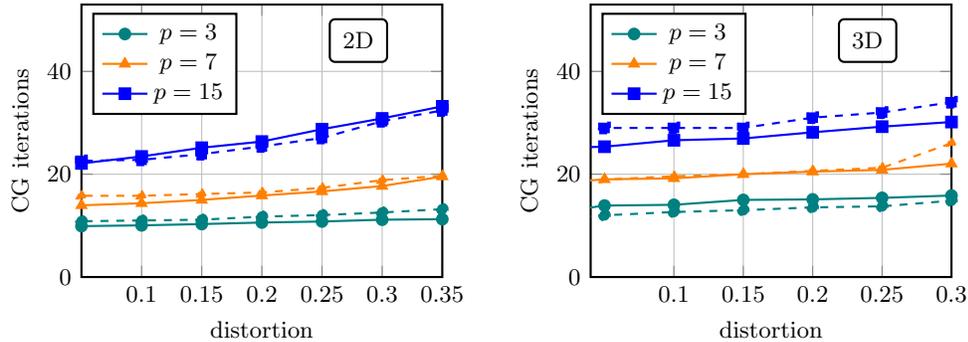
\begin{figure}[htbp]
  \centering
  \begin{tikzpicture}
    \begin{scope}
      \begin{axis}[
          paperplot,
          xlabel={distortion},
          ylabel={CG iterations},
          xmin=0.05,
          xmax=0.35,
          xtick={0.0,0.1,0.15,0.2,0.25, 0.3,0.35,0.4},
          ymin=0, ymax=53,
          legend pos=north west,
        ]
        \node[anchor=north east,inner sep=5pt,fill=white,fill opacity=0.9,text opacity=1,draw=black,rounded
          corners=2pt] at (rel axis cs:0.85,0.95) {2D};
        \pgfplotsset {cycle list set=0}
        \addplot+[ thick] table[
            col sep=comma,
            x=distortion,
            y=avg,
            restrict expr to domain={\thisrow{degree}}{2.9:3.1},
          ] {results/pureJac_ref_2D.csv};
        \addlegendentry{$p=3$}

        \addplot+[ thick] table[
            col sep=comma,
            x=distortion,
            y=avg,
            restrict expr to domain={\thisrow{degree}}{6.9:7.1},
          ] {results/pureJac_ref_2D.csv};
        \addlegendentry{$p=7$}

        \addplot+[ thick] table[
            col sep=comma,
            x=distortion,
            y=avg,
            restrict expr to domain={\thisrow{degree}}{14.9:15.1},
          ] {results/pureJac_ref_2D.csv};
        \addlegendentry{$p=15$}

        \pgfplotsset {cycle list set=0}

        \addplot+[ thick, dashed] table[
            col sep=comma,
            x=distortion,
            y=avg,
            restrict expr to domain={\thisrow{degree}}{2.9:3.1},
            restrict expr to domain={\thisrow{mu}}{0.9:1.1},
            unbounded coords=discard
          ] {results/jump_simplex_2D.csv};

        \addplot+[ thick, dashed] table[
            col sep=comma,
            x=distortion,
            y=avg,
            restrict expr to domain={\thisrow{degree}}{6.9:7.1},
            restrict expr to domain={\thisrow{mu}}{0.9:1.1},
            unbounded coords=discard
          ] {results/jump_simplex_2D.csv};

        \addplot+[ thick, dashed] table[
            col sep=comma,
            x=distortion,
            y=avg,
            restrict expr to domain={\thisrow{degree}}{14.9:15.1},
            restrict expr to domain={\thisrow{mu}}{0.9:1.1},
            unbounded coords=discard,
          ] {results/jump_simplex_2D.csv};

      \end{axis}
    \end{scope}

    \begin{scope}[xshift=0.52\columnwidth]
      \begin{axis}[
          paperplot,
          xlabel={distortion},
          ylabel={CG iterations},
          xmin=0.04, xmax=0.3,
          xtick={0.0,0.1,0.15,0.2,0.25, 0.3,0.35,0.4},
          ymin=0, ymax=53,
          legend pos=north west,
        ]
        \node[anchor=north east,inner sep=5pt,fill=white,fill opacity=0.9,text opacity=1,draw=black,rounded
          corners=2pt] at (rel axis cs:0.85,0.95) {3D};

        \pgfplotsset {cycle list set=0}
        \addplot+[ thick] table[
            col sep=comma,
            x=distortion,
            y=avg,
            restrict expr to domain={\thisrow{degree}}{2.9:3.1},
          ] {results/pureJac_simplex_3D.csv};

        \addplot+[ thick] table[
            col sep=comma,
            x=distortion,
            y=avg,
            restrict expr to domain={\thisrow{degree}}{6.9:7.1},
          ] {results/pureJac_simplex_3D.csv};

        \addplot+[ thick] table[
            col sep=comma,
            x=distortion,
            y=avg,
            restrict expr to domain={\thisrow{degree}}{14.9:15.1},
          ] {results/pureJac_simplex_3D.csv};

        \pgfplotsset {cycle list set=0}
        \addplot+[ thick, dashed] table[
            col sep=comma,
            x=distortion,
            y=avg,
            restrict expr to domain={\thisrow{degree}}{2.9:3.1},
          ] {results/simplex_Jac_cart_3D.csv};
        \addlegendentry{$p=3$}

        \addplot+[ thick, dashed] table[
            col sep=comma,
            x=distortion,
            y=avg,
            restrict expr to domain={\thisrow{degree}}{6.9:7.1},
          ] {results/simplex_Jac_cart_3D.csv};
        \addlegendentry{$p=7$}

        \addplot+[ thick, dashed] table[
            col sep=comma,
            x=distortion,
            y=avg,
            restrict expr to domain={\thisrow{degree}}{14.9:15.1},
          ] {results/simplex_Jac_cart_3D.csv};
        \addlegendentry{$p=15$}

      \end{axis}
    \end{scope}

  \end{tikzpicture}
  \caption{Number of CG iterations versus mesh distortion for the 2D(left) and 3D(right) patch solver when using a
    (damped) Jacobi smoother
    inside the p-multigrid V-cycle. Curves correspond to polynomial degrees $p=3,7,15$; each point shows the average
    number of CG iterations required to reduce the relative residual to $10^{-8}$. Solid lines show results for
    Cartesian patches; dashed lines show results for non-structured (simplex) patches in 2D.}
  \label{fig:iter_vs_distortion_pureJac}
\end{figure}

For a direct comparison, Figure~\ref{fig:iter_vs_distortion_pureJac} shows the performance of the simpler Jacobi
smoother on distorted patches. The left panel compares structured patches (solid lines) and unstructured simplex patches
(dashed lines) in 2D; the right panel shows results for 3D simplex patches.  In contrast to the Cartesian-reinforced
smoother, the (damped) Jacobi smoother is geometry‑agnostic and therefore more flexible: it remains applicable to arbitrary and
unstructured patches and shows only a gradual increase in iteration counts with distortion.  By comparing the results with the ones from Figure~\ref{fig:iter_vs_distortion}, it is evident that the Cartesian-reinforced smoother is
clearly superior on structured or near-Cartesian patches, delivering fewer iterations when its assumptions hold.
This demonstrates the practical trade-off: Cartesian-reinforced smoothers give best performance on well-structured grids,
while Jacobi smoothers offer flexibility as they do not rely on specific mesh structure.

\paragraph{Resilience to Variable Coefficients}
We now investigate robustness to variations in the material coefficient, focusing on a challenging discontinuous jump.
In this test, one cell in the patch is assigned a coefficient $\mu$ while all other cells have coefficient 1. Numerical
experiments indicate that the Cartesian-reinforced smoother, although effective for near-structured patches with mild
coefficient variation, fails to converge for contrast ratios exceeding roughly $\approx 10^{2}$. Consequently, for
these high-contrast cases, we employ the damped Jacobi smoother inside the p-multigrid cycle, which remains stable at
the cost of increased iteration counts.

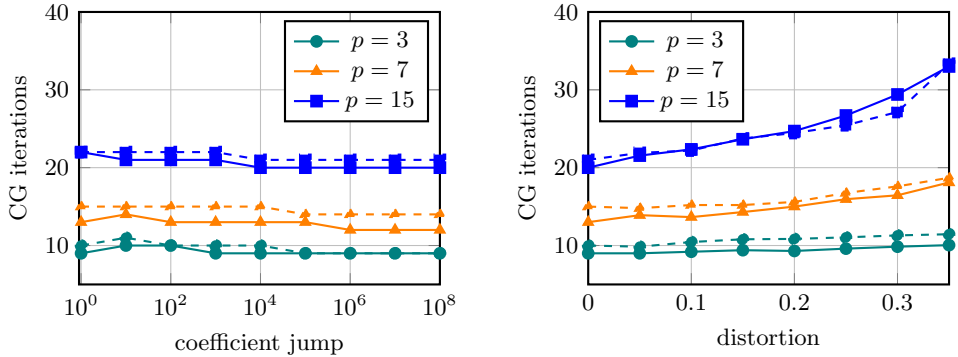
\begin{figure}[htbp]
  \centering
  \begin{tikzpicture}
    \begin{scope}
      \begin{axis}[
          paperplot,
          xlabel={coefficient jump},
          ylabel={CG iterations},
          xmin=.9,
          xmax=1e8,
          xmode=log,
          log basis x=10,
          ymin=5, ymax=40,
          legend pos=north east,
        ]
        \pgfplotsset {cycle list set=0}
        \addplot+[thick] table[
            col sep=comma,
            x=mu,
            y=avg,
            restrict expr to domain={\thisrow{degree}}{2.9:3.1},
            restrict expr to domain={\thisrow{distortion}}{-0.1:0.001},
          ] {results/jump_jac_2D.csv};
        \addlegendentry{$p=3$}

        \addplot+[ thick] table[
            col sep=comma,
            x=mu,
            y=avg,
            restrict expr to domain={\thisrow{degree}}{6.9:7.1},
            restrict expr to domain={\thisrow{distortion}}{-0.1:0.01},
          ] {results/jump_jac_2D.csv};
        \addlegendentry{$p=7$}

        \addplot+[thick] table[
            col sep=comma,
            x=mu,
            y=avg,
            restrict expr to domain={\thisrow{degree}}{14.9:15.1},
            restrict expr to domain={\thisrow{distortion}}{-0.1:0.01},
          ] {results/jump_jac_2D.csv};
        \addlegendentry{$p=15$}

        \pgfplotsset {cycle list set=0}

        \addplot+[ thick, dashed] table[
            col sep=comma,
            x=mu,
            y=avg,
            restrict expr to domain={\thisrow{degree}}{2.9:3.1},
            restrict expr to domain={\thisrow{distortion}}{-0.1:0.001},
          ] {results/jump_simplex_2D.csv};

        \addplot+[ thick, dashed] table[
            col sep=comma,
            x=mu,
            y=avg,
            restrict expr to domain={\thisrow{degree}}{6.9:7.1},
            restrict expr to domain={\thisrow{distortion}}{-0.1:0.01},
          ] {results/jump_simplex_2D.csv};

        \addplot+[ thick, dashed] table[
            col sep=comma,
            x=mu,
            y=avg,
            restrict expr to domain={\thisrow{degree}}{14.9:15.1},
            restrict expr to domain={\thisrow{distortion}}{-0.1:0.01},
          ] {results/jump_simplex_2D.csv};

      \end{axis}
    \end{scope}

    \begin{scope}[xshift=0.52\columnwidth]
      \begin{axis}[
          paperplot,
          xlabel={distortion},
          ylabel={CG iterations},
          xmin=0,
          xmax=0.35,
          ymin=5, ymax=40,
          legend pos=north west,
        ]

        \pgfplotsset {cycle list set=0}
        \addplot+[ thick, solid] table[
            col sep=comma,
            x=distortion,
            y=avg,
            restrict expr to domain={\thisrow{degree}}{2.9:3.1},
            restrict expr to domain={\thisrow{mu}}{9999:10001},
            unbounded coords=discard,
          ] {results/jump_jac_2D.csv};
        \addlegendentry{$p=3$}

        \addplot+[ thick] table[
            col sep=comma,
            x=distortion,
            y=avg,
            restrict expr to domain={\thisrow{degree}}{6.9:7.1},
            restrict expr to domain={\thisrow{mu}}{9999:10001},
            unbounded coords=discard,
          ] {results/jump_jac_2D.csv};
        \addlegendentry{$p=7$}

        \addplot+[ thick] table[
            col sep=comma,
            x=distortion,
            y=avg,
            restrict expr to domain={\thisrow{degree}}{14.9:15.1},
            restrict expr to domain={\thisrow{mu}}{9999:10001},
            unbounded coords=discard,
          ] {results/jump_jac_2D.csv};
        \addlegendentry{$p=15$}

        \pgfplotsset {cycle list set=0}
        \addplot+[ thick, dashed] table[
            col sep=comma,
            x=distortion,
            y=avg,
            restrict expr to domain={\thisrow{degree}}{2.9:3.1},
            restrict expr to domain={\thisrow{mu}}{9999:10001},
            unbounded coords=discard,
          ] {results/jump_simplex_2D.csv};

        \addplot+[ thick, dashed] table[
            col sep=comma,
            x=distortion,
            y=avg,
            restrict expr to domain={\thisrow{degree}}{6.9:7.1},
            restrict expr to domain={\thisrow{mu}}{9999:10001},
            unbounded coords=discard,
          ] {results/jump_simplex_2D.csv};

        \addplot+[ thick, dashed] table[
            col sep=comma,
            x=distortion,
            y=avg,
            restrict expr to domain={\thisrow{degree}}{14.9:15.1},
            restrict expr to domain={\thisrow{mu}}{9999:10001},
            unbounded coords=discard,
          ] {results/jump_simplex_2D.csv};

      \end{axis}
    \end{scope}
  \end{tikzpicture}
  \caption{
    Left panel: Number of CG iterations vs.\ coefficient jump on an undeformed 2D patch. All cells except one have
    coefficient $1$, while the remaining cell has coefficient $\mu$. Solid lines are for Cartesian patches, dashed for
    simplex patches.
    Right panel: Number of CG iterations vs.\ mesh distortion for a fixed coefficient jump of $\mu=10^4$ in 2D.
    In both plots, local patch problems are solved with a single p-multigrid V-cycle using a damped Jacobi smoother.
  }
  \label{fig:iter_vs_jump_pureJac}
\end{figure}

Figure~\ref{fig:iter_vs_jump_pureJac} presents two related experiments in 2D. The left panel plots the number of CG
iterations against the magnitude of the coefficient jump ($\mu$) on an undeformed patch. The results are very
encouraging: the number of iterations shows only a very weak dependence on the jump magnitude, growing very slowly even
for jumps spanning several orders of magnitude (from $10^0$ to $10^8$). This demonstrates that the p-multigrid
preconditioner with a Jacobi smoother is highly robust with respect to strong coefficient heterogeneity. The solid and
dashed lines show results for Cartesian and simplex patches, respectively, confirming the performance is consistent
across different element types.

The right panel examines the most challenging case, where the patch is both geometrically distorted and features a
large coefficient jump. It shows the results for a patch with a fixed, large coefficient jump ($\mu=10^4$) under
increasing geometric distortion. The iteration counts are close to those in the constant-coefficient case
(Figure~\ref{fig:iter_vs_distortion_pureJac}) and remain remarkably stable as the distortion increases. This indicates
that the solver's resilience to geometric distortion is not significantly compromised by the presence of a large
coefficient jump. Again, the trends are consistent for both Cartesian (solid lines) and unstructured simplex (dashed
lines) patches.

\paragraph{Summary}
The p-multigrid V-cycle with a simple Jacobi smoother proves to be a versatile and reliable choice. It handles
unstructured (simplex)
patches, large coefficient jumps, and geometric distortions with a graceful degradation in performance. This makes it a
great choice for patches with strong coefficient heterogeneity. However, its performance with respect to the polynomial
degree $p$ is suboptimal, with iteration counts roughly doubling as $p$ doubles. While already very good, this
suggests that there is room for improvement, and we are curious to see how this behavior translates to its performance
as a smoother in a full geometric multigrid solver.

On the other hand, the Cartesian-reinforced smoother is clearly superior for constant-coefficient problems on
near-structured grids, delivering excellent performance with very few iterations. Its effectiveness, however,
diminishes with increasing geometric distortion and it fails to converge for even moderate coefficient jumps. This
suggests that while powerful, its applicability is limited. A more advanced basis, such as the one proposed in
\cite{brubeck2021scalable}, could potentially serve as a better foundation for a specialized smoother that is more
resilient to such challenges. Overall, the smoothers are not fully robust to distortion, but they are resilient enough
that we expect them to perform well within a full multigrid hierarchy, which we will investigate next.

\section{Application as a smoother in geometric multigrid}
\label{sec:application_in_mg}

Having established the p-multigrid method as a robust solver for the local patch problem, we now return to its
application as a smoother within the global geometric multigrid algorithm. In a practical multigrid setting, the local
problem on each patch does not need to be solved to high accuracy. In fact, performing too much work in the smoother
can be inefficient, as it may not significantly improve the overall convergence. The goal of the smoother is to damp
high-frequency error
components, which can typically be achieved with just a few iterations of an effective iterative method.

Therefore, for the remainder of this paper, we denote by  $N_{\text{MG}}$ the number of p-multigrid V-cycles applied
within each local patch solve. Since the local solver serves as a smoother rather than an exact solver, we restrict
$N_{\text{MG}}$ to small values, specifically $N_{\text{MG}} \in \{1, 2\}$. This choice balances computational cost
with smoothing effectiveness, treating the local solver as an approximate inverse. The subsequent sections will analyze
the performance of the overall multigrid method when this approximate p-multigrid patch smoother is employed.
As before, in all experiments, we choose the damping parameter $\omega = \nicefrac{1}{2}$.

The test problems are run on a unit square (2D) or cube (3D), initially discretized by a single cell. This mesh  is then
uniformly refined to create a multigrid hierarchy. For the experiments on distorted Cartesian grids
(Tables~\ref{tab:gmres_iterations_jacobi} and~\ref{tab:gmres_iterations_cartesian_reinforced}), we use a fixed grid with
five levels in 2D and three levels in 3D. The corresponding problem sizes are highlighted in Table~\ref{tab:dof_counts_cartesian}.

\begin{figure}[htbp]
  \centering
  \begin{subfigure}[t]{0.32\textwidth}
    \centering
    \includegraphics[width=\textwidth]{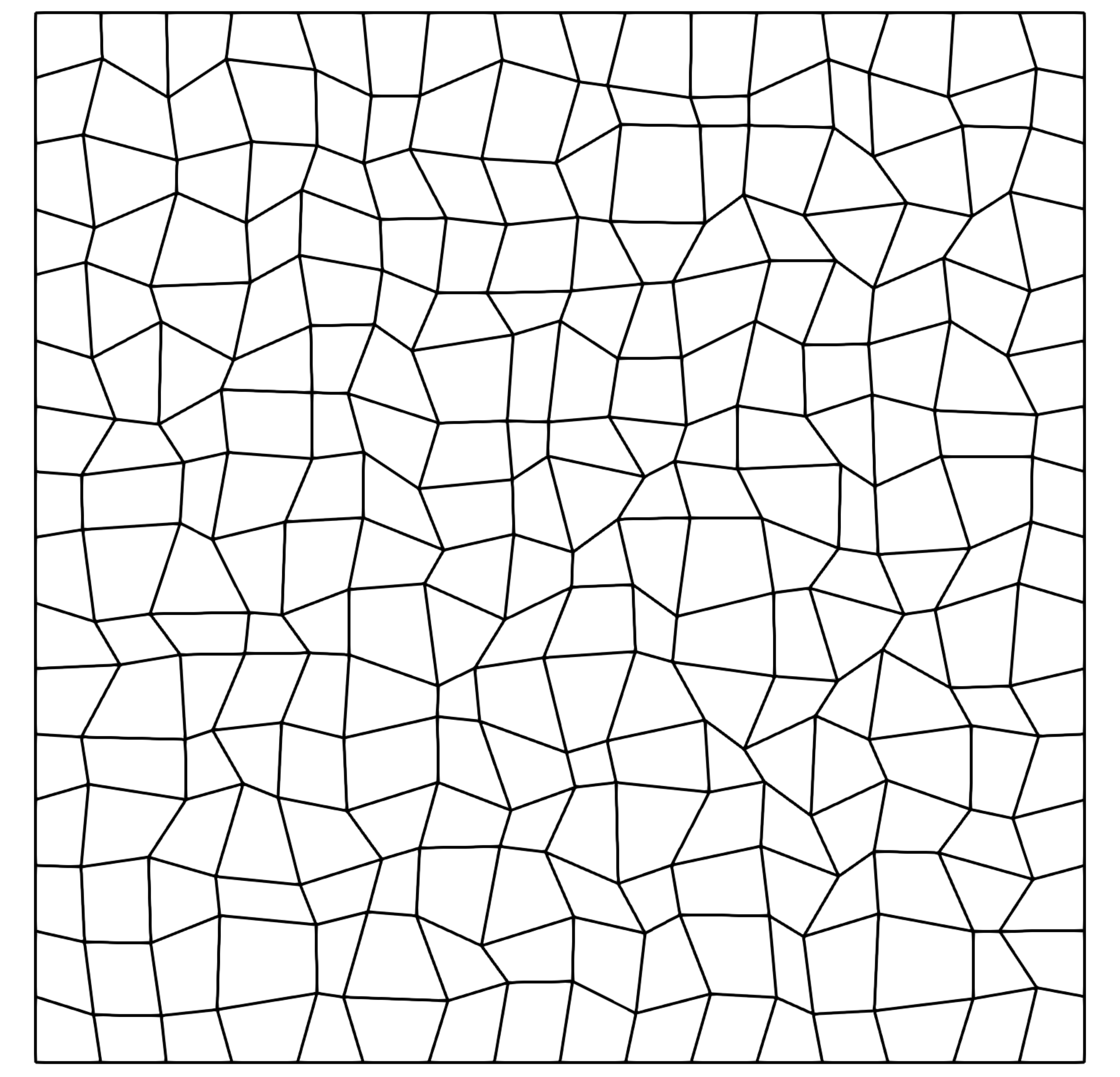}
    \caption{Mesh (2D)  distorted by $\delta= 30\%$ with  four global levels and}
    \label{fig:distorted_30}
  \end{subfigure}
  \hfill
  \begin{subfigure}[t]{0.32\textwidth}
    \centering
    \includegraphics[width=\textwidth]{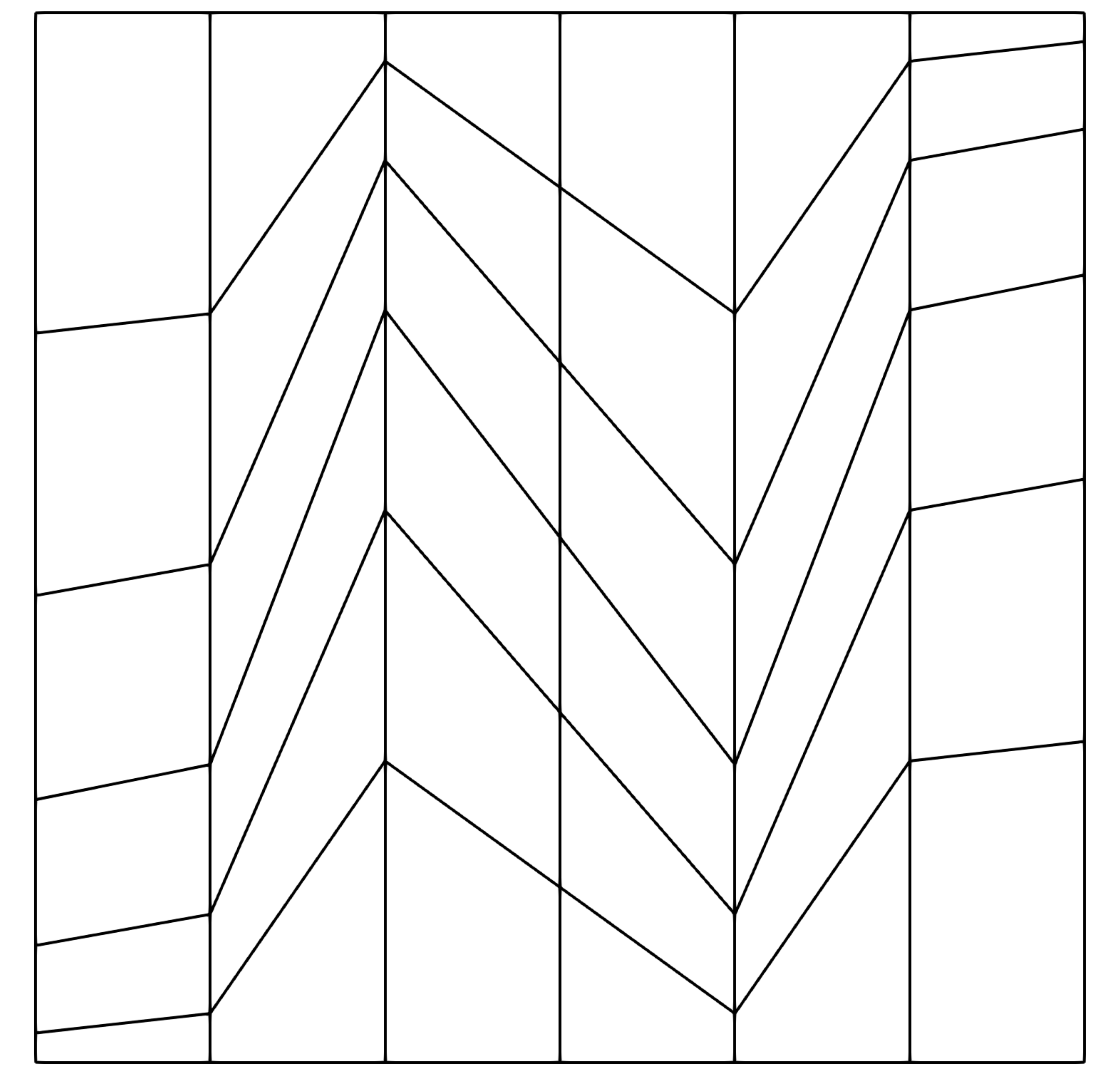}
    \caption{Kershaw 2D coarse mesh.}
    \label{fig:pickles_2d}
  \end{subfigure}
  \hfill
  \begin{subfigure}[t]{0.32\textwidth}
    \centering
    \includegraphics[width=\textwidth]{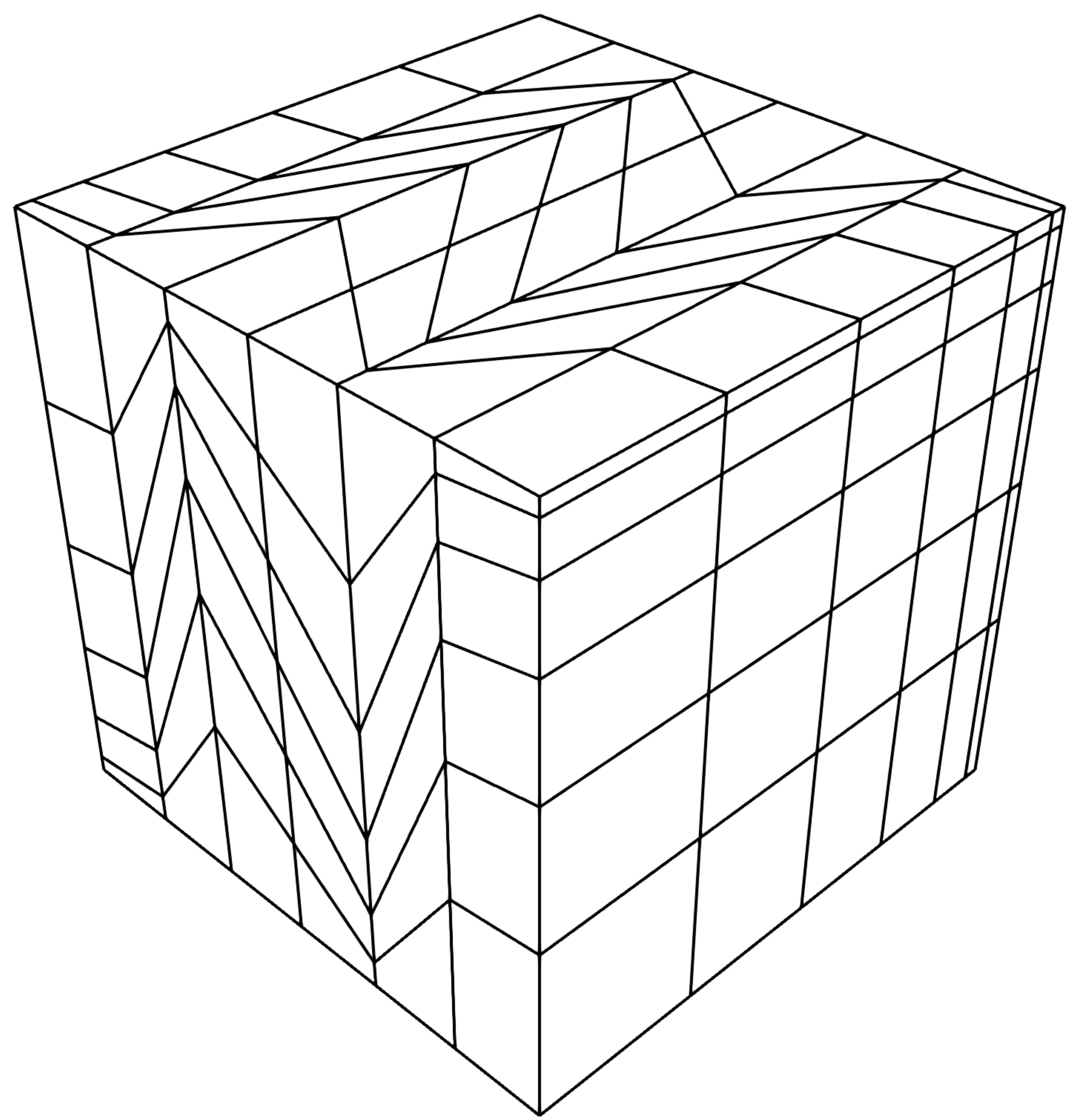}
    \caption{Kershaw 3D coarse mesh.}
    \label{fig:pickles_3d}
  \end{subfigure}
  \caption{Meshes used in the numerical experiments.}
\end{figure}

The test domain is constructed by starting with a single patch, refining it four times, and then distorting the
interior vertices. Each interior vertex $v$ is shifted by a distance $\delta h_v$ in a randomly chosen direction, where
$h_v$ is the minimum characteristic length of the edges attached to $v$ and $\delta$ is the distortion factor. A mesh
with $\delta=30\%$ is shown in Figure~\ref{fig:distorted_30}. Note that the mesh used in computations is refined
four times and then distorted in the same way. Some quadrilateral cells are distorted so heavily that they are nearly
triangular, with one interior angle close to 180 degrees.
In addition to these systematically distorted meshes, we also test our method on Kershaw-type grids~\cite{kershaw1981differencing}, which are known for their highly anisotropic and distorted cells, providing a challenging benchmark. The coarse meshes for these tests are shown in Figures~\ref{fig:pickles_2d} and~\ref{fig:pickles_3d}, and the corresponding problem sizes are given in Table~\ref{tab:dof_counts_kershaw}.
\definecolor{myhighlightcolor}{gray}{0.9}

\begin{table}[htbp]
  \centering
  \footnotesize
  \caption{Degrees of freedom (DoFs) for representative polynomial degrees $p$ on the finest grid level for Cartesian meshes. The highlighted columns correspond to the main test cases in 2D ( $L=5$) and 3D ($L=3$).}
  \label{tab:dof_counts_cartesian}
  \begin{tabular}{l|cccc|ccc}
    \toprule
    \multicolumn{1}{c}{} & \multicolumn{4}{c}{Cartesian 2D} & \multicolumn{3}{c}{Cartesian 3D}                                                                                                 \\
    $p\setminus L$       & 2                                & 3                                & 4     & \cellcolor{myhighlightcolor}5     & 2     & \cellcolor{myhighlightcolor}3     & 4     \\ 
    \midrule
    3                    & 169                              & 625                              & 2.40k & \cellcolor{myhighlightcolor}9.41k & 2.20k & \cellcolor{myhighlightcolor}15.6k & 118k  \\
    7                    & 841                              & 3.25k                            & 12.8k & \cellcolor{myhighlightcolor}50.6k & 24.4k & \cellcolor{myhighlightcolor}185k  & 1.44M \\
    15                   & 3.72k                            & 14.6k                            & 58.1k & \cellcolor{myhighlightcolor}231k  & 227k  & \cellcolor{myhighlightcolor}1.77M & 14.0M \\
    \bottomrule
  \end{tabular}
\end{table}

The patch smoother is applied multiplicatively, and the fixed patch ordering (which is the same for pre- and
post-smoothing) makes the global preconditioner non-symmetric. We therefore use GMRES as the outer iterative solver and
report the number of iterations required for convergence. The multigrid hierarchy always employs a single smoothing
step. The local patch solves use a p-multigrid method (i.e. without CG), so each local solver defines a linear
operator, even though the multiplicative composition makes the global preconditioner non-symmetric.
As the coarse solver we use GMRES preconditioned with a single sweep of the patch smoother. We solve the
system with a random right-hand side and report the number of iterations required to reduce the (relative) residual by
a factor of $10^{-8}$.

First, we consider the Jacobi smoother within the local p-multigrid solver. Table~\ref{tab:gmres_iterations_jacobi}
reports the GMRES iteration counts for various polynomial degrees $p$ and distortion factors. We test the performance
when applying one or two local p-multigrid V-cycles ($N_{\text{MG}}=1$ or $N_{\text{MG}}=2$) per patch. The results
show that even a single local V-cycle yields a robust smoother, with iteration counts remaining low across all tested
scenarios. For instance, in 2D with $p=15$ and 30\% distortion, the global solver requires only 22 iterations. This is
notably lower than the approximately $30$ iterations observed in the single-patch experiments under similar conditions (see
Figure~\ref{fig:iter_vs_distortion_pureJac}). We hypothesize that this improvement stems from the overlapping nature of
the patches: each cell is part of multiple patches, so its solution is updated several times within a single global
smoothing sweep, accelerating convergence.

To gauge the impact of the inexact local solve, we also report results for $N_{\text{MG}}=25$. This large number of
local cycles ensures that the local problem is solved to a high accuracy, making the local solver's inexactness
practically irrelevant. Comparing $N_{\text{MG}}=1$ or $2$ to $N_{\text{MG}}=25$, we see that a nearly exact local
solve consistently reduces the global iteration count, but the gains diminish as more local work is done. For example,
in 2D with $p=15$ and 30\% distortion, increasing $N_{\text{MG}}$ from 1 to 2 reduces iterations from 22 to 14, while a
further increase to 25 only brings it down to 10. This confirms that a small number of local V-cycles provides an
effective balance between computational cost and smoothing quality.

\begin{table}[htbp]
  \centering
  \centering
  \footnotesize
  \caption{GMRES iteration counts for a geometric multigrid preconditioner using a patch smoother with a local
    p-multigrid solver based on a Jacobi smoother. Each entry shows the number of iterations for an inexact local solve
    ($N_{\text{MG}}=1$ or $2$ p-MG V-cycles) versus a nearly exact local solve ($N_{\text{MG}}=25$ p-MG
    V-cycles).}
  \begin{tabular}{l c |cccc|cccc}
    \toprule
        &                                   & \multicolumn{4}{c|}{2D} & \multicolumn{4}{c}{3D}                                          \\
    $p$ & $N_{\text{MG}}  \setminus \delta$ & 0\%                     & 10\%                   & 25\% & 35\% & 0\% & 10\% & 25\% & 30\% \\
    \midrule
        & 1                                 & 5                       & 6                      & 7    & 9    & 6   & 6    & 7    & 7    \\
    3   & 2                                 & 4                       & 5                      & 6    & 7    & 4   & 5    & 5    & 5    \\
    \myrowcolor
        & 25                                & 4                       & 5                      & 6    & 7    & 4   & 4    & 5    & 5    \\
    \midrule
        & 1                                 & 5                       & 6                      & 8    & 10   & 8   & 8    & 9    & 9    \\
    7   & 2                                 & 4                       & 5                      & 6    & 8    & 6   & 6    & 7    & 7    \\
    \myrowcolor
        & 25                                & 3                       & 5                      & 6    & 7    & 3   & 4    & 5    & 5    \\
    \midrule
        & 1                                 & 6                       & 7                      & 9    & 11   & 10  & 10   & 11   & 12   \\
    15  & 2                                 & 5                       & 6                      & 7    & 9    & 7   & 8    & 8    & 9    \\
    \myrowcolor
        & 25                                & 3                       & 5                      & 6    & 7    & 3   & 4    & 5    & 5    \\
    \bottomrule
  \end{tabular}
  \label{tab:gmres_iterations_jacobi}
\end{table}

Next, we investigate the Cartesian-reinforced smoother, with results presented in
Table~\ref{tab:gmres_iterations_cartesian_reinforced}. This smoother demonstrates impressive resilience to the
polynomial degree $p$, even on highly distorted meshes. While there is a slight growth in iteration counts with $p$,
this trend is mirrored when using a nearly exact local solver (the values in parentheses, corresponding to
$N_{\text{MG}}=25$). This suggests that the moderate increase in iterations is an intrinsic property of the smoother on
distorted grids, rather than an artifact of the inexact local solve. As expected, the Cartesian-reinforced smoother
consistently
outperforms the Jacobi smoother on these structured-grid problems, confirming its superiority when its underlying
This result is particularly noteworthy, as it was questioned in~\cite{brubeck2021scalable} whether p-robustness could
be maintained for highly distorted meshes. Our findings demonstrate that this is indeed possible. For example, as shown
in Table~\ref{tab:gmres_iterations_cartesian_reinforced}, for a highly distorted mesh (30\%) in 3D, the number of GMRES
iterations increases only from 8 to 10 as the polynomial degree is raised from $p=3$ to $p=15$.

\begin{table}[htbp]
  \centering
  \centering
  \footnotesize
  \caption{GMRES iteration counts for a geometric multigrid preconditioner using a patch smoother with a local
    p-multigrid solver based on a Cartesian-reinforced smoother. Each entry shows the iteration count for an inexact
    local
    solve ($N_{\text{MG}}=1$ p-MG V-cycle) and, in parentheses, for a nearly exact local solve ($N_{\text{MG}}=25$ p-MG
    V-cycles). The bar --- indicates no convergence within 25 iterations.
  }
  \begin{center}
    \begin{tabular}{l|cccc|cccc}
      \toprule
                          & \multicolumn{4}{c|}{2D} & \multicolumn{4}{c}{3D}                                                   \\
      $p\setminus \delta$ & 0\%                     & 10\%                   & 25\%  & 35\%    & 0\%   & 10\%  & 25\%  & 30\%  \\
      \midrule
      3                   & 4 (4)                   & 5 (5)                  & 6 (6) & 8 (7)   & 5 (4) & 5 (4) & 6 (5) & 6 (5) \\
      \myrowcolor
      7                   & 3 (3)                   & 5 (5)                  & 6 (6) & 9 (7)   & 4 (3) & 4 (4) & 5 (5) & 6 (5) \\
      15                  & 3 (3)                   & 5 (5)                  & 6 (6) & --- (7) & 4 (3) & 4 (4) & 5 (5) & 6 (5) \\
      \bottomrule
    \end{tabular}
  \end{center}

  \label{tab:gmres_iterations_cartesian_reinforced}
\end{table}

To place our results in context, we compare them with other state-of-the-art methods. The work of Witte et
al.~\cite{WitteArndtKanschat21} provides a useful benchmark, employing a cell-wise additive smoother on meshes with
25\% distortion. For
a polynomial degree of $p=15$, they reported 66 iterations in 2D (five levels) and 80 in 3D (three levels). Under
similar conditions, our multiplicative vertex-patch smoother with the Cartesian-reinforced local solver is
significantly more efficient, as detailed in Table~\ref{tab:gmres_iterations_cartesian_reinforced}.

A direct comparison is also possible with the vertex-patch method of Brubeck et al.~\cite{brubeck2021scalable},
where an additive smoother on Kershaw grids was used. The results for our method on similar grids are
presented in Table~\ref{tab:hp-comparison}. On a highly distorted Kershaw grid with two multigrid levels,
they reported 158 iterations for $p=7$ in 3D . For a comparable setup (L=2), our method requires
15 iterations. This stark contrast highlights the efficiency of the proposed p-multigrid-based smoother.
Furthermore, Table~\ref{tab:hp-comparison} demonstrates that our method is not only p-robust but
also h-robust, as the iteration counts tend to decrease with an increasing number of multigrid levels (L).

\begin{table}[tbhp]
  {
    \footnotesize
    \caption{
      GMRES iteration counts for a geometric multigrid preconditioner using a patch smoother with a local
      p-multigrid solver based on a Jacobi smoother with $N_{\text{MG}}=1$. The results for Kershaw meshes with $p=15$ have not been computed due to high coarse solver cost. Numbers in parentheses are from~\cite{brubeck2021scalable} for comparison.
    }
    \label{tab:hp-comparison}
    \begin{center}
      \begin{tabular}{rl|ccc|ccc|ccc}
        \toprule
            &                                     & \multicolumn{3}{c|}{Distorted 10\%}
            & \multicolumn{3}{c|}{Distorted 25\%}
            & \multicolumn{3}{c}{Kershaw}                                                                                            \\
        $d$ & $p\setminus L$
            & 2                                   & 3                                   & 4
            & 2                                   & 3                                   & 4
            & 2                                   & 3                                   & 4                                          \\
        \midrule
            & 7                                   & 8                                   & 7  & 7 & 9  & 9  & 8  & 16 (54)  & 15 & 13 \\ \myrowcolor
        2   & 3                                   & 6                                   & 6  & 6 & 7  & 7  & 7  & 13 (78)  & 14 & 14 \\
            & 15                                  & 11                                  & 10 & 8 & 12 & 11 & 10 & 19 (90)  & 15 & 13 \\
        \midrule
            & 7                                   & 9                                   & 8  & 7 & 10 & 9  & 8  & 18 (102) & 18 & 16 \\ \myrowcolor
        3   & 3                                   & 7                                   & 6  & 6 & 7  & 7  & 7  & 15 (158) & 17 & 17 \\
            & 15                                  & 12                                  & 10 & 9 & 13 & 11 & 10 &          &    &    \\
        \bottomrule
      \end{tabular}
    \end{center}
  }
\end{table}

\begin{table}[htbp]
  \centering
  \footnotesize
  \caption{Degrees of freedom (DoFs) for representative polynomial degrees $p$ on the finest grid level for Kershaw meshes.}
  \label{tab:dof_counts_kershaw}
  \begin{tabular}{l|ccc|ccc}
    \toprule
    \multicolumn{1}{c}{} & \multicolumn{3}{c}{Kershaw 2D} & \multicolumn{3}{c}{Kershaw 3D}                                 \\
    $p\setminus L$       & 2                              & 3                              & 4     & 2     & 3     & 4     \\
    \midrule
    3                    & 1.36k                          & 5.32k                          & 21.0k & 50.7k & 389k  & 3.05M \\
    7                    & 7.22k                          & 28.5k                          & 113k  & 614k  & 4.83M & 38.3M \\
    15                   & 32.7k                          & 130k                           & 519k  &       &       &       \\
    \bottomrule
  \end{tabular}
\end{table}

\section{Conclusion}
We have explored matrix-free patch smoothers for geometric multigrid methods; a powerful approach that
has attracted considerable attention in recent years. Like many before us, we encountered the fundamental challenge of
efficiently solving the local problems that arise on each patch. Rather than resorting to direct factorization or other
specialized techniques, we followed the recursive philosophy inherent to multigrid methods and
thus constructed the local solvers using another multigrid solver, employing p-coarsening this time. This recursive
\emph{multigrid-within-multigrid} approach yields a matrix-free local solver that remains true to the iterative
spirit of the method.

Our numerical experiments revealed distinct personalities among the smoothers. The p-multigrid V-cycle with Jacobi
smoothing proved to be a reliable generalist: it handled unstructured patches, large coefficient jumps, and geometric
distortions with consistent performance. The Cartesian-reinforced smoother, by contrast, showed its specialist
nature: delivering exceptional results on structured grids: for a mesh with 25\% distortion in 3D, the global GMRES iteration count drop from 6 to 5 as the polynomial degree increases from $3$ to $15$.

The Jacobi-based smoother provided a particularly encouraging answer on the Kershaw grids that were found to be particularly challenging in~\cite{brubeck2021scalable}, significantly outperforming results reported there. In~\cite{brubeck2021scalable} also questioned whether p-robustness could survive highly distorted meshes. Our results with the Jacobi-based smoother show it can: for a Kershaw grid in 3D with 4 multigrid levels, the global GMRES iteration count grows only from 16 to 17 as the polynomial degree increases from $3$ to $7$. We have thus
demonstrated a smoother that achieves robustness on heavily distorted meshes while maintaining optimal local
computational complexity of $O(p^{d+1})$ and memory requirements of $O(p^{d})$.

Beyond theoretical interest, the method offers practical advantages. It requires minimal precomputation, achieves
optimal memory usage, and adapts to both structured and unstructured meshes. By following the recursion one, or
perhaps a couple levels deeper, we arrived at a smoother that is robust, efficient, and well-suited for high-order
finite element applications.

\section*{Acknowledgments}
I would like to thank Guido Kanschat for the  motivation he provided for this work.

The generative AI  (Gemini, ChatGPT, Claude) was used to assist in drafting and proofreading
parts of this manuscript. Any remaining errors are the author's responsibility.


\FloatBarrier

\bibliographystyle{siamplain}
\bibliography{literature}
\end{document}